\documentclass[12pt]{amsart}
\usepackage[dvips]{epsfig}
\usepackage{graphics}
\usepackage{latexsym}
\usepackage{verbatim}
\usepackage{amsmath}
\usepackage{amsthm}
\usepackage{amssymb}
\newtheorem{theorem}{Theorem}[section]
\newtheorem{lemma}[theorem]{Lemma}

\newtheorem{proposition}[theorem]{Proposition}

\newtheorem{remark}[theorem]{Remark}

\def\Remark{\medskip\noindent{\bf Remark: }}

\newcommand{\ens}[1]{\mathbb{#1}}

\newcommand{\N}{\mathbb{N}}

\newcommand{\R}{\mathbb{R}}

\def\derpar#1#2{\frac{\partial#1}{\partial#2}}
\def\var{\varepsilon}
\def\pa{\partial}

\begin{document}

\title[Stability for the Boltzmann equation 
with long-range interactions]
{Stability and uniqueness for the spatially homogeneous 
Boltzmann equation with long-range interactions}

\author{Laurent Desvillettes, Cl\'ement Mouhot}

\hyphenation{bounda-ry rea-so-na-ble be-ha-vior pro-per-ties
cha-rac-te-ris-tic}

\begin{abstract} 
In this paper, we prove some {\it a priori} stability estimates (in weighted Sobolev spaces) 
for the spatially homogeneous Boltzmann equation without angular cutoff (covering every 
physical collision kernels). These estimates are conditioned to some regularity estimates on 
the solutions, and therefore reduce the stability and uniqueness issue to the one of proving suitable
regularity bounds on the solutions. We then prove such regularity bounds for a class of interactions
including the so-called (non cutoff and non mollified) hard potentials and moderately soft potentials.
In particular, we obtain the first result of global existence and uniqueness for these long-range 
interactions.
\end{abstract}

\maketitle



\tableofcontents

\section{Introduction}
\setcounter{equation}{0}

\subsection{The Boltzmann equation}

The Boltzmann equation (Cf. \cite{Ce88} and \cite{CIP}) 
describes the behavior of a dilute gas when the only 
interactions taken into account are binary collisions.
In the case when the distribution function is
assumed to be independent on the position $x$, we obtain 
the so-called {\it spatially homogeneous Boltzmann equation}, which reads
 \begin{equation}\label{el}
 \derpar{f}{t}(t,v)  = Q(f,f)(t,v), \qquad  v \in \R^N, \quad t \geq 0,
 \end{equation}
where $N \ge 2$ is the dimension. 
In equation~\eqref{el}, $Q$ is the quadratic
{\it Boltzmann collision operator}, defined by the bilinear symmetrized form
 \begin{equation*}\label{eq:collop}
 Q(g,f)(v) = \frac12\,\int _{\R^N \times \ens{S}^{N-1}} B(|v-v_*|, \cos \theta)
         \left(g'_* f' + g' f_* '- g_* f - g f_* \right) \, dv_* \, d\sigma,
 \end{equation*}
where we have used the shorthands $f=f(v)$, $f'=f(v')$, $g_*=g(v_*)$ and
$g'_*=g(v'_*)$. Moreover, $v'$ and $v'_*$ are parametrized by
 \begin{equation*}\label{eq:rel:vit}
 v' = \frac{v+v_*}2 + \frac{|v-v_*|}2 \, \sigma, \qquad 
 v'_* = \frac{v+v_*}2 - \frac{|v-v_*|}2 \, \sigma, \qquad \sigma \in \ens{S}^{N-1}. 
 \end{equation*}
Finally, $\theta\in [0,\pi]$ is the deviation angle between 
$v'-v'_*$ and $v-v_*$ defined by $\cos \theta = (v'-v'_*)\cdot(v-v_*)/|v-v_*|^2$, 
and $B$ is the Boltzmann collision kernel determined by physics 
(related to the cross-section $\Sigma(v-v_*,\sigma)$ 
by the formula $B=|v-v_*| \, \Sigma$).  We also formally denote  
$$ Q^+(f,f)(v) =  \int _{\R^N \times \ens{S}^{N-1}} B(|v-v_*|, \cos \theta)\, f'_* f'\, dv_* \, d\sigma$$
the gain part of $Q$, and 
$$ L(f)(v) = \int _{\R^N \times \ens{S}^{N-1}} B(|v-v_*|, \cos \theta)\, f_*\, dv_* \, d\sigma $$
the linear operator appearing in the loss part $Q^-$ of $Q$.
\medskip

Boltzmann's collision operator has the fundamental properties of
conserving mass, momentum and energy
  \begin{equation}
  \int_{\R^N}Q(f,f) \, \phi(v)\,dv = 0, \qquad
  \phi(v)=1,v,|v|^2, \label{CON}
  \end{equation}
and satisfying Boltzmann's $H$ theorem, which writes (at the formal level) 
  \begin{equation*} 
  - \frac{d}{dt} \int_{\R^N} f \log f \, dv = - \int_{\R^N} Q(f,f)\log(f) \, dv \geq 0.
  \end{equation*}

\subsection{Assumptions on the collision kernel}  

We shall consider the following assumptions on the collision kernel $B$: 

 \begin{itemize}
 \item[{\bf H1.}] It takes the following tensorial form (with $\Phi,b$ nonnegative functions)
	\[ B(|v-v_{*}|, \cos \theta) = \Phi(|v-v_{*}|) \, b(\cos \theta). \]

\smallskip
 \item[{\bf H2.}] The angular part is nonnegative, smooth (or at least locally integrable) 
 for $\theta \in (0,\pi]$, and such that
    \[ b(\cos \theta) \sim_{\theta \to 0} C_b \, \theta^{-(N-1)-\nu} \]
  with $\nu \in (-\infty,2)$ and $C_b>0$. 
 \end{itemize}

As for the ``kinetic'' part $\Phi$, we make one of the following assumptions:
\begin{itemize}
 \item[{\bf H3-1.}] The function  $z \mapsto \Phi(|z|)$ is strictly positive, 
 $C^{\infty}$, such that 
 $$ \Phi(|z|) \sim_{|z| \to +\infty} C_\Phi \, |z|^{\gamma}, $$
 for some $C_{\Phi} > 0$ and $\gamma \in (0,1]$, and satisfies the bounds 
   \[  \forall \, z \in \R^N, \ p\in \N^*,  \quad |\Phi^{(p)} (|z|) | \le C_{\Phi,p}, \]
 for some $C_{\Phi,p} >0$. 

\smallskip
 \item[{\bf H3-2.}] The function  $z \mapsto \Phi(|z|)$ is strictly positive, 
 $C^{\infty}$, such that 
 $$ \Phi(|z|) \sim_{|z| \to +\infty} C_\Phi \, |z|^{\gamma} $$
 for some  $C_{\Phi}>0$ and $\gamma \in (-N,0]$, and satisfies the bounds 
   \[  \forall \, z \in \R^N, \ p\in \N^*,  \quad |\Phi^{(p)} (|z|) | \le C_{\Phi,p}, \]
 for some $C_{\Phi,p} >0$. 

\smallskip
 \item[{\bf H3-3.}] The function  $z \mapsto \Phi(|z|)$ is given by the explicit formula
       \[\Phi(|z|) =  C_\Phi \, |z|^{\gamma}, \]
 for some  $C_{\Phi}>0$ and $\gamma \in (0,1]$.
 
\smallskip
 \item[{\bf H3-4.}] The function  $z \mapsto \Phi(|z|)$ is given by the explicit formula
       \[\Phi(|z|) =  C_\Phi \, |z|^{\gamma}, \]
 for some  $C_{\Phi}>0$ and $\gamma \in (-N,0]$.
 \end{itemize}

Our assumptions (more precisely, {\bf H3-3}) cover in dimension $3$ the hard
spheres collision kernel $B(|v-v_*|, \cos \theta)= \mbox{{\rm cst}} \, |v-v_*|$.
It also covers (still in dimension $3$) collision kernels deriving from interaction 
potentials behaving like inverse-power laws. More precisely 
for an interaction potential $V(r) = \mbox{{\rm cst}} \, r^{-s}$, $B$ satisfies our 
assumptions with the formulas $\gamma = (s-4)/s$ and $\nu = 2/s$ (see~\cite{Ce88}). 
One traditionally calls {\em hard potentials} the case $s >4$ 
(for which $0< \gamma <1$, and which corresponds to {\bf H3-3}),
{\em Maxwell molecules} the case $s=4$ (which corresponds to {\bf H3-2} with $\gamma=0$),
and {\em soft potentials} the case $1<s<4$ (for which $-N<\gamma<0$, and which
corresponds to {\bf H3-4}). 
\par
Assumptions {\bf H3-1} and {\bf H3-2} correspond to cases when $B$ is artificially smoothed around
$0$ with respect to $v-v_*$. 
\par
Since $\nu = 2/s$ for potentials in $r^{-s}$, only the non-negative $\nu$ are physically
meaningful (as far as inverse power laws are concerned). The case of negative $\nu$, 
corresponding to the so-called {\em angular cutoff}, is a simplification.

\subsection{Goals, existing results and difficulties} 

The stability of the spatially homogeneous  Boltzmann equation for hard potentials (or hard spheres,
or Maxwellian molecules) {\em with angular cutoff} was proven, in weighted $L^1$ spaces, by Arkeryd~\cite{Ark72}. The special structure of the Maxwellian molecules makes it possible to prove the
stability (and consequently the uniqueness) of the corresponding spatially homogeneous  Boltzmann equation
{\em without angular cutoff} using Fourier transform and Wasserstein-like distances (Cf. \cite {toscanivillaniwasserstein}). 
A recent work of Fournier~\cite{Fournier} shows by probabilistic means that this stability also 
holds {\em without angular cutoff} for ``kinetic'' sections $\Phi$  which are not of Maxwellian molecules type but
 are bounded and smooth (this more or less corresponds to our hypothesis {\bf H3-2}),
 and for moderate angular singularities (that is $\nu \in ]0,1[$).
 This recent paper is an important step in the application of contraction metrics approach for non constant collision kernels,
 even if up to now it does not consider physical collision kernels apart from the Maxwell molecules one.

Hence, as far as we know, no stability (or uniqueness) result is known for ``true'' hard or soft potentials.  
Here, we show that stability holds for any kind of interactions, as soon as suitable regularity bounds are at hand.
 Then, we prove these required regularity bounds for a class of models including ``true'' hard potentials and
 moderately soft potentials.   
Our approach is complementary to the one of Fournier~\cite{Fournier} in the sense that our stability
is in a stronger space (that is in a weighted $W^{1,1}$ space instead of a measure space), for
more general cross-sections, 
but leads to a uniqueness result holding for a smaller set of initial data. 
Our method of proof is also completely different.
It is based on the use of integrations by parts for finite differences of a special kind.

We shall devote a separate forthcoming work~\cite{DMasymp} to the question of the asymptotic behavior when $t\to +\infty$
 of the solution of the spatially homogeneous Boltzmann equation {\em without cutoff},
 on the basis of the Cauchy theorems established in the present work, new entropy production estimates,
 and the approach developed in~\cite{DM:05}. 

\subsection{Notation} 

We denote $\langle \cdot \rangle = (1 + |\cdot|^2)^{1/2}$.  
We shall systematically use the following notations ($s\in\R$, $p\in [1, +\infty)$, $k \in \N$)
\[ \|f\|_{L^p_s}^p := \int_{\R^N} |f(v)|^p\, \langle v \rangle^{ps} \, dv, \qquad
\|f\|_{L^{\infty}_s} := \sup_{v\in \R^N} |f(v)|\, \langle v \rangle^{s}\]
and 
\[ \|f\|^p _{W^{k,p} _s} := \sum_{0\le |i| \le k} \|\partial^i f\|_{L^p _s}^p,  \qquad
 \|f\|_{W^{k,\infty} _s} := \sum_{0\le |i| \le k} \|\partial^i f\|_{L^{\infty} _s}, \]
where $\partial^i$ denotes the partial derivative related to the 
multi-index $i$. In the case $p=2$, we denote $H^k _s = W^{k,2} _s$. 
We finally use the notation $x_+$ for the nonnegative part of $x \in \R$, defined by 
$x_+ = \max \{ x, 0 \}$. 

\subsection{Statement of the results}

We first state the key {\em a priori}  stability theorem for moderate angular singularities:
 
 \begin{theorem}\label{theo:stab}
 Let $B$ be a collision kernel which satisfies {\bf H1-H2} with $\nu <1$, and let  $f,g \in L^{\infty}([0,T]; L^1 _2 \cap L\log L (\R^N))$ be two nonnegative solutions to the spatially homogeneous Boltzmann equation 
 associated to $B$, on some time interval $[0,T]$.

We assume first that $B$ satisfies {\bf H3-1}, {\bf H3-2} or {\bf H3-3}. For any $q \ge 2$,
we have the following {\em a priori} bound 
 \begin{equation} \label{GGH}
 \forall \, t \in [0,T], \quad \|f(t,\cdot) - g(t,\cdot)\|_{L^1 _q} \le 
          \|f_0 - g_0\|_{L^1 _q} \, \exp( C_s \, t ), 
	  \end{equation}
 with 
   \[ C_s = \mbox{{\rm cst}} \, \left( \sup_{t \in [0,T]} \| f(t,\cdot)\|_{W^{1,1} _{q+(1+\gamma)_+}} 
      + \sup_{t \in [0,T]} \| g(t,\cdot)\|_{W^{1,1} _{q+ (1+\gamma)_+}} \right) .  \]

We assume then that $B$ satisfies {\bf H3-4}. We still have (\ref{GGH}) for any $q \ge 2$, but with
   \[ C_s = \mbox{{\rm cst}} \, \left( \sup_{t \in [0,T]} \| f(t,\cdot)\|_{W^{1,1}_{q+(1 + \gamma)_+} \cap L^{p}} 
  + \sup_{t \in [0,T]} \| g(t,\cdot)\|_{W^{1,1} _{q+(1+ \gamma)_+} \cap L^{p}} \right) , \]
where $p > N/(N+\gamma)$ if $\gamma+1\ge 0$, and else 
 \[ C_s = \mbox{{\rm cst}} \, \left( \sup_{t \in [0,T]} \| f(t,\cdot)\|_{W^{1,1} _{q+(1 + \gamma)_+} \cap L^{p_1}} 
     + \sup_{t \in [0,T]} \| \nabla f (t,\cdot)\|_{L^{p_2}} \right.   \]
  \[ \left. + \sup_{t \in [0,T]} \| g(t,\cdot)\|_{W^{1,1} _{q+(1+ \gamma)_+} \cap L^{p_1}}
 + \sup_{t \in [0,T]}  \| \nabla g(t,\cdot) \|_{L^{p_2}}  \right) , \]
where $p_1 > N/(N+\gamma)$ and $p_2 > N/(N+\gamma+1)$. 
\end{theorem}
\bigskip

We also give a proposition stating the stability result for strong angular singularities. 
For the sake of simplicity, we do not write down the explicit estimate in this case for 
``true'' soft potentials (case {\bf H3-4}), but it can be obtained from our proof.  
\bigskip

\begin{proposition}\label{propo:sta}
Let $B$ be a collision kernel which satisfies {\bf H1-H2} with $ 1 \le \nu < 2$, and
let $0 \le f(t, \cdot), g(t,\cdot) \in L^1 _2$ be two solutions to the spatially homogeneous Boltzmann equation
(\ref{el}) associated to $B$, on some time interval $[0,T]$. We assume first that $B$
 satisfies {\bf H3-1}, {\bf H3-2} or {\bf H3-3}. For any $q \ge 4$, estimate (\ref{GGH}) holds with
   \[ C_s = \mbox{{\rm cst}} \, \left( \sup_{t \in [0,T]} \| f(t,\cdot)\|_{W^{2,1} _{q+(2+\gamma)_+}} 
      + \sup_{t \in [0,T]} \| g(t,\cdot)\|_{W^{2,1} _{q+ (2+\gamma)_+}} \right) . \]
\end{proposition}

Second, we state a theorem summing up what we obtain for the Cauchy theory by combining the previous {\em a priori} 
 stability estimates (that is, Theorem~\ref{theo:stab})
 with results on the propagation of smoothness  which are either already known 
 (when $\nu < 0$, that is, for cutoff cross-sections)  or new 
(when $\nu \in [0,1)$, that is for non cutoff cross-sections).

\begin{theorem}\label{theo:reg}
Let $B$ be a collision kernel which satisfies {\bf H1-H2} with  $\nu <1$. 
\begin{itemize} 
\item If $B$ satisfies {\bf H3-1}, {\bf H3-2} or {\bf H3-3} and $f(0,\cdot)$ is an initial datum belonging to 
$W^{1,1} _q$ for some $q \ge 2$, there is a unique global solution to eq. (\ref{el}) in the space 
$W^{1,1} _q$. 
\item If $B$ satisfies {\bf H3-4} for $\gamma \ge -1$, and $f(0,\cdot)$ is an initial datum belonging to
$W^{1,1} _q \cap L^p$ for $p > \frac{N}{N+\gamma}$ and $q\ge 2$,
  there is a unique local (that is, on a certain time interval $[0,T]$)
 solution in the space $W^{1,1} _q \cap L^p$. Moreover, this solution is global (that is, $T= +\infty$) 
when $\gamma \in (-\nu,0]$ and $q$ is big enough (depending on $\gamma$, $\nu$).
\end{itemize}
\end{theorem}

\begin{remark} \label{rem}
Hence in dimension $N=3$ where $\nu = 2/s$ and $\gamma = (s-4)/s$ for potentials in $r^{-s}$, this theorem yields 
global existence and uniqueness results for any ``true'' hard potentials ($4< s < +\infty$) and ``true'' 
moderately soft potentials ($2 < s < 4$). Indeed in these cases, one has $0< \nu <1$ and $-\nu< \gamma \le 1$. 
Note however that our method, like that of Fournier, does not seem to work for strong angular singularities
($1<\nu < 2$), even if the kinetic part $\Phi$ of the cross-section is very smooth. As a consequence, 
Proposition~\ref{propo:sta} has not yet found an application (that is,
it gives a result of uniqueness without existence\dots )
\end{remark}


\subsection{Plan of the paper}

Section \ref{pse} is devoted to the proof of Theorem~\ref{theo:stab} and
Proposition~\ref{propo:sta}. The case of ``true'' soft potentials {\bf H3-4}
is somewhat different from the others, and it is therefore treated separately.
\par
Then, Theorem~\ref{theo:reg} is proven in section \ref{eps}. Once again, the
case  of ``true'' soft potentials deserves a special treatment.

\section{Proof of the stability estimates}\label{pse}
\setcounter{equation}{0}

{\bf{Proof of  Theorem~\ref{theo:stab} and Proposition~\ref{propo:sta}}:}
Without restriction (since $Q$ is taken in symmetrized form), we replace 
in the whole paper the cross-section $B$ by its symmetrized form with 
support included in $\theta \in [0,\pi/2]$:
 \begin{equation*} 
 B_{\mbox{{\scriptsize sym}}} (|v-v_*|, \cos \theta) = \Big[ B(|v-v_*|, \cos\theta) + 
 B(|v-v_*|, \cos(\pi-\theta)) \Big] \, {\bf 1}_{\cos\theta\geq 0},
 \end{equation*}
where ${\bf 1}_E$ denotes the usual characteristic function of the set $E$. 
\par
Let $B$ be a collision kernel which satisfies {\bf H1-H2} and one of the assumptions
{\bf H3}, and 
$f(t,\cdot)$, $g(t,\cdot) \ge 0$ be two solutions on $[0,T]$ to the associated spatially homogeneous Boltzmann 
equation (\ref{el}).
Let us define  
$D=f-g$  and $S=f+g$. 
The evolution equation for $D$ reads 
 \[ \derpar{D}{t} = Q(f,f) - Q(g,g) = Q(S,D). \]

\subsection{Hard potentials or mollified soft potentials} 
\par
Let us first assume {\bf H3-1}, {\bf H3-2} or {\bf H3-3}.

We split $b = b_c ^{\var} + b_r ^{\var}$ with 
$b_c ^\var  = b \, {\bf 1}_{\theta \in [\var,\pi/2]}$ (the cutoff part), 
$b_r ^\var= 1-b_c ^\var$ (the remaining part), and $\var \in (0,\pi/2]$ to be fixed later. 
This induces corresponding splittings of the collision kernel $B = B_c + B_r$ and of the 
collision operator $Q = Q_c + Q_r$.
 
Then, we have
 \begin{multline*}
 \frac{d}{dt} \|D\|_{L^1 _q} \le 
   \int_{\R^N} Q_r(S,D) \, \mbox{sgn}(D) \, \langle v \rangle^q \, dv 
    + \int_{\R^N}  
   Q_c(S,D) \, \mbox{sgn}(D) \, \langle v \rangle^q \, dv 
                             =: I_1 + I_2.
 \end{multline*}

For the cutoff part,  we estimate
 \begin{multline*}
 2\,I_2 = 2\,\int_{\R^N} Q_c(S,D) \, \mbox{sgn}(D) \, \langle v \rangle^q \, dv \\
 = \int_{v,v_*,\sigma} \left[ S' _* D' + S' D' _* - S_* D - S D_* \right] 
 \, \mbox{sgn}(D) \, \langle v \rangle^q \, B_c \\
 \le \int_{v,v_*,\sigma} \left( S' _* |D'| + S' |D'| _* - S_* |D| - S |D_*| \right) 
 \, \langle v \rangle^q \, B_c
 + 2\,\int_{v,v_*,\sigma} S |D_*| \, \langle v \rangle^q \, B_c \\
 =  2 \,\int_{\R^N}  Q_c(S,|D|) \, \langle v \rangle^q \, dv 
 + 2 \,\int_{v,v_*,\sigma} S |D_*| \, \langle v \rangle^q \, B_c,       
 \end{multline*}
which implies 
 \begin{multline*}
 I_2 \le  C\,\int_{\R^N}  Q_c(S,|D|) \, \langle v \rangle^q \, dv \\
       + C_\var \, \int_{\R^N \times \R^N} |D_*| S \, 
                    \langle v \rangle^{q+\gamma_+} \, \langle v_* \rangle^{\gamma_+} \, dv_* \, dv 
                   =: I_{2,1} + I_{2,2}, 
 \end{multline*}
where the constant $C_\var >0$ depends on $\var >0$ {\em via} the $L^1$ norm of $b_c ^\var$ 
on the sphere $\ens{S}^{N-1}$ (which possibly blows up as $\var \to 0$).

The $I_{2,2}$ term is controlled (for $q \ge \gamma_+$) by 
 \[ I_{2,2} \le C _\var \, C_q \, \|D\|_{L^1 _{\gamma_+}} 
 \le C' _\var \, \|D\|_{L^1 _q}. \]

The $I_{2,1}$ term writes, using the pre-post-collisional change 
of variable (see~\cite[Chapter~1, Section~4.5]{Vill:hand}):
 \[ I_{2,1} = \int_{\R^N \times \R^N} |D_*| S \, \Phi(|v-v_*|) \,
 \left( \int_{\ens{S}^{N-1}} 
                 \big[ \langle v' \rangle^q
         + \langle v' _* \rangle^q - \langle v \rangle^q - \langle v_* \rangle^q \big] 
                \, b_c ^\var \, d\sigma \right) \, dv \, dv_*. \]
Then, we shall prove (for $q >2$)
 \begin{multline*}
 \left( \int_{\ens{S}^{N-1}} 
                 \Big[ \langle v' \rangle^q
           + \langle v' _* \rangle^q - \langle v \rangle^q - \langle v_* \rangle^q \Big] 
                \, b_c ^\var \, d\sigma \right) \\  
 \le  \mbox{{\rm cst}} \,\left( \int_{\ens{S}^{N-1}} 
                 \Big[ |v'|^q + |v' _*|^q - |v|^q - |v_*|^q \Big] 
		 \, b_c ^\var \, d\sigma \right) \\
 \le C_\var \, \Big[ |v|^{q-1} \, |v_*| + |v| \, |v_*|^{q-1} \Big] 
          - K \Big[ |v|^q + |v_*|^q \Big] 
 \end{multline*}
for again some constant $C_\var >0$ possibly blowing-up as $\var \to 0$, and some  
constant $K>0$ independent of $\var$. 

The proof of this inequality is straightforward, using the kind of Povzner inequalities 
in~\cite{Wenn:momt:97,LuWe:02}. Indeed, \cite[Lemma~1]{LuWe:02} implies that 
\begin{multline*}
\langle v' \rangle^q + \langle v' _* \rangle^q - \langle v \rangle^q - \langle v_* \rangle^q \\
           \le 2^{q+1} \, \Big[ \langle v \rangle^{q-1} \, \langle v_* \rangle 
                  + \langle v_* \rangle^{q-1} \, \langle v \rangle \Big] \, \cos \theta \, \sin \theta 
                - K_q \, \Big[ \langle v \rangle^q  + \langle v_* \rangle^q \Big] \, \cos^2 \theta \, \sin^2 \theta 
\end{multline*}
for some constant $K_q>0$ depending only on $q$ (note that the proof in~\cite{LuWe:02}  is done in dimension $3$
 but straightforwardly extends to any dimension). Note also that 
in the case of moderate angular singularities $\nu < 1$, the constant $C_\var$ indeed does not blow up as $\var$ goes to infinity. 
 
Hence, using that $q \ge 1 + \gamma_+$, that $q-1 + \gamma_+ \le q$ (since $\gamma_+ \le 1$ 
and $q>2$), and also that 
$$
\int_{\R^N} S(v_*) \, \Phi(v-v_*) \, dv_* \ge \mbox{{\rm cst}} \, \langle v \rangle^{\gamma} 
$$
thanks 
to the entropy bounds on $f$ and $g$, 
we get
 \[  I_{2,1} \le C^{'}_\var \, C_q \, \|D\|_{L^1 _q} - K' \, \|D\|_{L^1 _{q+\gamma}}, \]
with $C' _\var$ possibly blowing up as $\var \to 0$,
 and $K'$ independent of $\var$. 
\bigskip

The remaining non cutoff part writes
 \begin{multline*}
 \int_{\R^N} Q_r(S,D) \, \mbox{sgn}(D) \, \langle v \rangle^q \, dv \\ 
    = \int_{\R^N \times \R^N \times \ens{S}^{N-1}} \Phi \, b_r \, 
       \Big[ S' _* D' - S_* D \Big]  \, \mbox{sgn}(D) \, \langle v \rangle^q \, dv \, dv_* \, d\sigma \\ 
    + \int_{\R^N \times \R^N \times \ens{S}^{N-1}} \Phi \, b_r \, 
       \Big[ D' _* S' - D_* S \Big]  \, \mbox{sgn}(D) \, \langle v \rangle^q \, dv \, dv_* \, d\sigma 
    =: I_{1,1} + I_{1,2}.
 \end{multline*}

The $I_{1,1}$ term is the easiest to deal with:
 \begin{multline*}
 I_{1,1} = \int_{v,v_*,\sigma} \big[ S' _* D' - S_* D \big] 
                  \, \mbox{sgn}(D) \, \langle v \rangle^q 
               \, \Phi(|v-v_*|) \, b_r ^\var \, dv \, dv_* \, d\sigma \\
     =  \int_{v,v_*,\sigma} S_* \Big[ D \, \mbox{sgn}(D') \langle v' \rangle^q 
                                       - |D| \langle v \rangle^q \Big] 
               \, \Phi(|v-v_*|) \, b_r ^\var \, dv \, dv_* \, d\sigma \\
    \le \int_{v,v_*,\sigma} S_* |D| \left[ \langle v' \rangle^q 
                                       - \langle v \rangle^q \right] 
               \, \Phi(|v-v_*|) \, b_r ^\var \, dv \, dv_* \, d\sigma \\ 
    \le \int_{v,v_*} S_* |D| \, \Phi(|v-v_*| )
               \left| \int_{\ens{S}^{N-1}} \Big[ \langle v' \rangle^q 
                                       - \langle v \rangle^q \Big] \,  
                \, b_r ^\var \, d\sigma \right| \, dv \, dv_*. 
  \end{multline*}

Then we shall prove a simple lemma, which is a variant of \cite[Lemma~2.3]{DM:03}. 

\begin{lemma}\label{diffpoids}
Let $q \ge 2$, then 
\begin{multline*}
\left| 
\int_{\ens{S}^{N-1}} \Big[ \langle v' \rangle^q - \langle v \rangle^q \Big] 
\, b(\cos \theta) \, d\sigma \right| \\
\le C \, \left( \int_{\ens{S}^{N-1}} b(\cos \theta) \, \sin \theta/2 \, d\sigma \right) \, 
|v - v_*| \, \Big[  \langle v \rangle^{q-1} + \langle v_* \rangle^{q-1} \Big]. 
\end{multline*}
Let $q \ge 4$, then 
\begin{multline*}
\left| \int_{\ens{S}^{N-1}} \Big[ \langle v' \rangle^q - \langle v \rangle^q \Big] 
\, b(\cos \theta) \, d\sigma \right| \\ 
\le C \, \left( \int_{\ens{S}^{N-1}} b(\cos \theta) \, \left( \sin \theta/2 \right)^2 \, d\sigma \right) \, 
|v - v_*|^2 \, \Big[  \langle v \rangle^{q-2} + \langle v_* \rangle^{q-2} \Big]. 
\end{multline*}
In those formulas, the constants $C >0$ depend only on $q$. The same formulas 
are true when $(v,v')$ is replaced by $(v_*,v' _*)$.
\end{lemma}  

\smallskip\noindent{\sl Proof of Lemma \ref{diffpoids}.} 
The proof is straightforward by using integral Taylor expansions of 
$\ u \in [0,1] \mapsto \langle v'_u \rangle^q \in \R$ 
(denoting $v' _u = (1- u) \, v + u \, v'$). 

At first order, one gets 
$$ 
\langle v' \rangle^q - \langle v \rangle^q = 
q \, \int_0 ^1 \langle v' _u \rangle^{q-2} \, v' _u \cdot (v'-v) \, du, 
$$
which is enough to prove the first inequality. 
\par
At second order, one gets
\begin{multline*}
\langle v' \rangle^q - \langle v \rangle^q = q \, \langle v \rangle^{q-2} \, v \cdot (v'-v) \\
+ q \, \int_0 ^1 \Big[ \langle v' _u \rangle^{q-2} \, |v'-v|^2 
+ (q-2) \,  \langle v' _u \rangle^{q-4} \, |v'_u|^2 \, |v'-v|^2 \Big] \, (1-u) \, du.  
\end{multline*}
This is enough to prove the second inequality as soon as one notices that 
$$
\int_{\ens{S}^{N-1}} \langle v \rangle^{q-2} \, v \cdot (v'-v) \, b(\cos \theta) \, d\sigma = 0
$$
by gathering antipodal points of the $(N-2)$-dimensional sphere $\ens{S}^{N-1} \cap (v-v_*)^\bot$. 
\qed
\bigskip

We turn back to the proof of Theorem~\ref{theo:stab} and Proposition~\ref{propo:sta}
in the case of hard or smoothed soft potentials.
\par
From Lemma~\ref{diffpoids},
we deduce that when $q\ge 2$,
 \[  \left| \int_{\ens{S}^{N-1}} \Big[ \langle v' \rangle^q 
                                       - \langle v \rangle^q \Big] \, b_r ^\var \, d\sigma \right| 
      \le C \, m_1(b_r ^\var) \, \langle v \rangle^q \, \langle v_* \rangle^q, \]
where $C >0$ is independent of $\var$, and
 \[ m_1(b_r ^\var) = \int_{\ens{S}^{N-1}} b_r ^\var (\cos \theta) \,\sin \theta/2 \, d\sigma \]
is a finite quantity which goes to $0$ as $\var$ goes to $0$,  
under the assumption $\nu \in (0,1)$. 
For $q \ge 4$, one also has the control 
 \[  \left| \int_{\ens{S}^{N-1}}  \Big[ \langle v' \rangle^q 
                                       - \langle v \rangle^q \Big]  \, b_r ^\var \, d\sigma \right| 
     \le C \, m_2(b_r ^\var) \, \langle v \rangle^q \, \langle v_* \rangle^q \]
with 
 \[ m_2(b_r ^\var) = \int_{\ens{S}^{N-1}} b_r ^\var (\cos \theta) 
     \, \left( \sin \theta/2 \right) ^2 \, d\sigma \]
which is finite and goes to $0$ as $\var$ goes to $0$, for any $\nu <2$ 
(that is the whole physical range).

Therefore assuming $q \ge 2$ when $\nu <1$ 
(case of Theorem \ref{theo:stab}) or $q \ge 4$ when 
$1\le \nu <2$ (case of Proposition~\ref{propo:sta}), 
we obtain that $I_{1,1}$ is controlled by $\|D\|_{L^1 _{q+\gamma_+}}$ times 
some constant which goes to $0$ as $\var$ goes to $0$. Thus 
 \[ I_{1,1} \le \frac{K'}{4} \, \|D\|_{L^1 _{q+\gamma_+}} \] 
for $\var$ small enough.

We now come to the most difficult term to estimate, and the crucial point in the proof. 
That is the use of suitable changes of variables which play (loosely speaking) the role 
of some integration by parts for the ``integral differentiation-like operators'' appearing 
in the collision operator for grazing collisions.  
\par
At this point, in order to keep tractable notations, we keep on with the proof
only under assumption {\bf H3-3}. It can be checked easily that the proof
also works under assumptions {\bf H3-1} and {\bf H3-2} (replacing $\gamma$ by
$\gamma_+$ if necessary, and using bounds on the derivatives of $\Phi$).
\par
 The $I_{1,2}$ term writes
 \begin{multline*}
 I_{1,2} = \int_{v,v_*,\sigma} \big[ D' _* S' - D_* S \big] 
                  \, \mbox{sgn}(D) \, \langle v \rangle^q 
               \, |v-v_*|^\gamma \, b_r ^\var \, dv \, dv_* \, d\sigma \\
     =  \int_{v,v_*,\sigma} D_* S \big[ \mbox{sgn}(D') \langle v' \rangle^q 
                                       - \mbox{sgn}(D) \langle v \rangle^q \big] 
               \,  |v-v_*|^\gamma \, b_r ^\var \, dv \, dv_* \, d\sigma \\
    \le \int_{v,v_*,\sigma} D_* \mbox{sgn}(D) \langle v \rangle^q 
                     \left( \frac{S(\phi_\sigma(v,v_*))}{(\cos \theta/2)^{N+\gamma}} - S \right)
               \,  |v-v_*|^\gamma \, b_r ^\var \, dv \, dv_* \, d\sigma \\
     \le \int_{v,v_*} |D_*| \langle v \rangle^q 
               \left| \int_{\ens{S}^{N-1}} \left( 
          \frac{S(\phi_\sigma(v,v_*))}{(\cos \theta/2)^{N+\gamma}} - S \right) 
              \, b_r ^\var \, d\sigma \right| \, |v-v_*|^\gamma \, dv \, dv_* ,
  \end{multline*}
where we have used the change of variable from cancellation lemmas 
in~\cite{ADVW} (which is possible since $b$ has its support included 
in $[0,\pi/2]$). The variable $\phi_\sigma(v,v_*)$ denotes the inverse application of 
$v \mapsto v'$ keeping $v_*$ and $\sigma$ frozen (it is given explicitly in~\cite{ADVW}).
Let us denote $\bar{v} = \phi_\sigma(v,v_*)$. 

We split the integral on the sphere into three parts
 \begin{multline*}
 \langle v \rangle^q \, \int_{\ens{S}^{N-1}}  \bigg(
           \frac{S(\bar{v})}{(\cos \theta/2)^{N+\gamma}} - S 
            \, \bigg) \, b_r ^\var \, d\sigma  \\
  = \int_{\ens{S}^{N-1}} 
     \frac{S(\bar{v}) \langle \bar{v} \rangle^q 
     - S(v) \langle v \rangle^q }{(\cos \theta/2)^{N+\gamma}}
     \, b_r ^\var \, d\sigma  \\
 + \int_{\ens{S}^{N-1}} 
     \left( \frac{1}{(\cos \theta/2)^{N+\gamma}} -1 \right) \, b_r ^\var \, d\sigma
     \, S(v) \, \langle v \rangle^q 
 + \int_{\ens{S}^{N-1}} 
     \left( \frac{\langle v \rangle^q - 
     \langle \bar{v} \rangle^q}{(\cos \theta/2)^{N+\gamma}} \right) 
     S(\bar{v}) \, b_r ^\var \, d\sigma ,
 \end{multline*}
which yields a corresponding splitting of $I_{1,2}$ into three parts 
$I_{1,2,1} + I_{1,2,2} + I_{1,2,3}$. 

For the $I_{1,2,3}$ term, we use again the change of variable defined above, but backward: 
 \[ |I_{1,2,3}| \le C \, \int_{v,v_*} |D_*| S(v) \, |v-v_*|^\gamma \,  
               \left( \int_{\ens{S}^{N-1}} 
          \left| \langle v' \rangle^q - \langle v \rangle^q \right| 
              \, b_r ^\var \, d\sigma \right) \, dv \, dv_* , \]
and we apply again Lemma~\ref{diffpoids} (with $q \ge 2$ if $\nu <1$ or $q \ge 4$ if $1 \le \nu <2$) 
to get  
by choosing $\var$ small enough
 \[  I_{1,2,3} \le \frac{K'}{4} \, \|D\|_{L^1 _{q+\gamma}}. \]

For the $I_{1,2,2}$ term, we have 
 \[ \left| \frac{1}{(\cos \theta/2)^{N+\gamma}} -1 \right| \le C \, (1-\cos \theta), \]
and so 
 \[  I_{1,2,2} \le C \, C_s \, m_2(b_r ^\var) \, \|D\|_{L^1 _{q+\gamma}}. \] 
Hence, choosing again $\var$ small enough, we get
 \[  I_{1,2,2} \le \frac{K'}{4} \, \|D\|_{L^1 _{q+\gamma}}. \]

Finally, for the $I_{1,2,1}$ term, we denote $\bar{v}_u = (1-u)\,v + u\,\bar{v}$ for 
$u \in [0,1]$ and we Taylor-expand $v \mapsto S(v) \langle v \rangle^q$. 
Let us first suppose that $\nu <1$. Then, it is enough to go to first order:  
$$
S(\bar v) \langle \bar v \rangle^q -S(v) \langle v \rangle^q 
= \int_0 ^1 \nabla(S \langle \cdot \rangle^q) (\bar{v}_u) \cdot (\bar v - v) \, du,
$$ 
and using the identity $|\bar v - v| = \tan \theta/2 \, |v - v_*|$, we get 
 \[ I_{1,2,1} \le C \, \int_{\R^{2N} \times \ens{S}^{N-1} \times [0,1]}
       |\nabla(S \langle \cdot \rangle^q)|(\bar{v}_u) \, 
       |v-v_*|^{\gamma+1} \, (\tan \theta/2) \, b_r ^\var
       \, |D_*| \, dv \, dv_* \, d\sigma \, du. \]
When $u$, $v_*$ and $\sigma$ are fixed, the change of variable 
$v \to \bar{v}_u$ has its Jacobian determinant bounded by a constant, 
and for any $u \in [0,1]$, 
 \[ |v-v_*|^{\gamma+1} \le C \, |\bar{v}_u-v_*|^{\gamma+1}, \]
hence
 \begin{multline*}
 I_{1,2,1} \le C \, \int_{[0,1]} \, \int_{\R^{2N} \times \ens{S}^{N-1}}
       |\nabla(S \langle \cdot \rangle^q)(v)| \, 
       |v-v_*|^{\gamma+1} \, (\tan \theta/2) \, b_r ^\var
       \, |D_*| \, dv \, dv_* \, d\sigma \, du \\
      \le C \, \left(\int_{\ens{S}^{N-1}}(\tan \theta/2) \, b_r ^\var \, d\sigma \right) \, 
           \|S\|_{W^{1,1} _{q+\gamma+1}} \, \|D\|_{L^1 _{\gamma +1}}, 
 \end{multline*}
for some constant $C>0$ independent of $\var$. 

Thus, remembering that $\nu <1$, we have 
 \[ \int_{\ens{S}^{N-1}}(\sin \theta/2) \, b_r ^\var \, d\sigma < +\infty\] 
and so (for $q \ge \gamma +1$)
 \[ I_{1,2,1} \le C \, C_s \, \|D\|_{L^1 _{q}}. \]
 
Let us now briefly explain  how to adapt this proof when $1 \le \nu <2$ (case of 
Proposition~\ref{propo:sta}).
 In order to cancel singularities 
of order $2$, we Taylor-expand at second order:
$$
S(\bar v) \langle \bar v \rangle^q -S(v) \langle v \rangle^q 
= \nabla(S \langle \cdot \rangle^q) (v) \cdot (\bar v - v)
+ \int_0 ^1  \nabla^2 (S \langle \cdot \rangle^q) (\bar{v}_u) \cdot (\bar v - v) \cdot (\bar v - v) \, du.
$$ 
Then, the key remark is that when $v,v_*,\theta$ are fixed, the unit vector $\sigma$ describes 
a sub-sphere of dimension $(N-2)$ included in $\ens{S}^{N-1}$, and the integral of the first order term 
over this sub-sphere is zero by gathering antipodal points. Therefore, we get 
 \[ I_{1,2,1} \le C \, \int_{\R^{2N} \times \ens{S}^{N-1} \times [0,1]}
       \big|\nabla^2(S \langle \cdot \rangle^q)\big|(\bar{v}_u) \, 
       |v-v_*|^{\gamma+2} \, (\tan \theta/2)^2 \, b_r ^\var
       \, |D_*| \, dv \, dv_* \, d\sigma \, du. \]
Using the same backward change of variable as in the case $\nu < 1$, we deduce 
 \begin{multline*}
 I_{1,2,1} \le C \, \int_{[0,1]} \, \int_{\R^{2N} \times \ens{S}^{N-1}}
       \big|\nabla^2(S \langle \cdot \rangle^q)(v)\big| \, 
       |v-v_*|^{\gamma+2} \, (\tan \theta/2)^2 \, b_r ^\var
       \, |D_*| \, dv \, dv_* \, d\sigma \, du \\
      \le C \, \left(\int_{\ens{S}^{N-1}}(\tan \theta/2)^2 \, b_r ^\var \, d\sigma \right) \, 
           \|S\|_{W^{2,1} _{q+\gamma+2}} \, \|D\|_{L^1 _{\gamma +2}}, 
 \end{multline*}
for some constant $C>0$ independent of $\var$.  Thus, since 
 \[ \int_{\ens{S}^{N-1}} (\sin \theta/2)^2 \, b_r ^\var \, d\sigma < +\infty,\] 
we get again (for $q \ge \gamma +2$)
 \[ I_{1,2,1} \le C \, C_s \, \|D\|_{L^1 _{q}}. \]

Combining all the previous estimates, we deduce (for $q \ge 2$  when $\nu < 1$ and $q\ge 4$ when $1 \le \nu < 2$)
 \[  \frac{d}{dt} \|D\|_{L^1 _q} \le C_+ \, \|D\|_{L^1 _q} - K_- \, \|D\|_{L^1 _{q+\gamma}}, \]
which concludes the proof of Theorem~\ref{theo:stab} and Proposition~\ref{propo:sta} in the case 
when assumptions {\bf{H3-1}}, {\bf{H3-2}} or {\bf{H3-3}} hold. 


\subsection{Non mollified soft potentials}
\par
We now assume that {\bf H3-4} holds (that is, in particular, $\gamma \le 0$).  We do not need to 
perform the splitting between cutoff and non-cutoff parts since large velocities are well-behaved,
but another difficulty occurs because of the singularity of the kinetic collision kernel $\Phi$ for 
small relative velocities.  We write the proof shortly, pointing out the differences with the previous subsection. 

We have again 
 \begin{multline*}
 2\, \int_{\R^N} Q(S,D) \, \mbox{sgn}(D) \, \langle v \rangle^q \, dv \\ 
    = \int_{\R^N \times \R^N \times \ens{S}^{N-1}} \Phi \, b \, 
       \big[ S' _* D' - S_* D \big]  \, \mbox{sgn}(D) \, \langle v \rangle^q \, dv \, dv_* \, d\sigma \\ 
    + \int_{\R^N \times \R^N \times \ens{S}^{N-1}} \Phi \, b \, 
       \big[ D' _* S' - D_* S \big]  \, \mbox{sgn}(D) \, \langle v \rangle^q \, dv \, dv_* \, d\sigma 
    =: I_{1} + I_{2}.
 \end{multline*} 
For the $I_{1}$ term again, we use Lemma~\ref{diffpoids} 
to deduce (assuming $q \ge 2$)
\begin{equation}\label{ll}
\bigg|\int_{\sigma} (\langle v' \rangle^q - \langle v \rangle^q ) \, b\, d\sigma \bigg|
\le \mbox{{\rm cst}} \, (\langle v \rangle^q + \langle v_* \rangle^q).
\end{equation}
Then we compute (for $p > N/(N+\gamma)$)
 \begin{multline*}
 I_{1} = \int_{v,v_*,\sigma} \big[ S' _* D' - S_* D \big] 
                  \, \mbox{sgn}(D) \, \langle v \rangle^q 
               \, |v-v_*|^\gamma \, b \, dv \, dv_* \, d\sigma \\
     =  \int_{v,v_*,\sigma} S_* \big[ D \, \mbox{sgn}(D') \langle v' \rangle^q 
                                       - |D| \langle v \rangle^q \big] 
               \, |v-v_*|^\gamma \, b \, dv \, dv_* \, d\sigma \\
    \le \int_{v,v_*,\sigma} S_* |D| \big[ \langle v' \rangle^q 
                                       - \langle v \rangle^q \big] 
               \, |v-v_*|^\gamma \, b \, dv \, dv_* \, d\sigma \\ 
    \le C \, \int_{v,v_*} S_* |D| \, |v-v_*|^\gamma \, 
        \big[ \langle v \rangle^q +  \langle v_* \rangle^q \big] \, dv \, dv_*
	\\
	\le C\, \int D\, \langle v \rangle^q \, \bigg( \int_{|v-v_*| \ge 1}  S_*  
	+ \int_{|v-v_*| \le 1}  S_*  \, |v-v_*|^\gamma  \bigg)  \\
	+  C\, \int D\,  \bigg( \int_{|v-v_*| \ge 1}  S_* \,  \langle v_* \rangle^q  
	+ \int_{|v-v_*| \le 1}  S_*  \, |v-v_*|^\gamma\, \langle v_* \rangle^q \bigg)  \\
	 \le C\, \int D\, \langle v \rangle^q \, \bigg( \|S\|_{L^1} + \| S \|_{L^p} \bigg)
	 +  C\, \int D\,  \bigg( \| S \|_{L^1_q} + \mbox{{\rm cst}} \, \langle v \rangle^q \, \| S \|_{L^p} \bigg)
	 \\
	 \le C \, \left( \|S\|_{L^1_q} + \|S\|_{L^p} \right) 
	                  \, \|D\|_{L^1 _q}. 
  \end{multline*}

The $I_{2}$ term writes again
 \begin{multline*}
 I_{2} = \int_{v,v_*,\sigma} \big[ D' _* S' - D_* S \big] 
                  \, \mbox{sgn}(D) \, \langle v \rangle^q 
               \, |v-v_*|^\gamma \, b \, dv \, dv_* \, d\sigma \\
     =  \int_{v,v_*,\sigma} D_* S \big[ \mbox{sgn}(D') \langle v' \rangle^q 
                                       - \mbox{sgn}(D) \langle v \rangle^q \big] \, 
              |v-v_*|^\gamma \, b \, dv \, dv_* \, d\sigma \\
    \le \int_{v,v_*,\sigma} D_* \mbox{sgn}(D) \langle v \rangle^q 
                     \left[ \frac{S(\phi_\sigma(v,v_*))}{(\cos \theta/2)^{N+\gamma}}
                                       - S \right] \, 
             |v-v_*|^\gamma \, b \, dv \, dv_* \, d\sigma \\
     \le \int_{v,v_*} |D_*| \langle v \rangle^q 
               \left| \int_{\ens{S}^{N-1}} 
          \bigg( \frac{S(\phi_\sigma(v,v_*))}{(\cos \theta/2)^{N+\gamma}} - S \bigg)
              \, b \, d\sigma \right| \,  |v-v_*|^\gamma \, dv \, dv_* ,
  \end{multline*}
where the spherical integral splits into 
 \begin{multline*}
    \langle v \rangle^q \,  \int_{\ens{S}^{N-1}} 
         \bigg(   \frac{S(\bar{v})}{(\cos \theta/2)^{N+\gamma}} - S \bigg)
            \, b \, d\sigma  \\
  = \int_{\ens{S}^{N-1}} 
     \frac{S(\bar{v}) \langle \bar{v} \rangle^q 
     - S(v) \langle v \rangle^q }{(\cos \theta/2)^{N+\gamma}}
     \, b \, d\sigma  \\
 + \int_{\ens{S}^{N-1}} 
     \left( \frac{1}{(\cos \theta/2)^{N+\gamma}} -1 \right) \, b \, d\sigma
     \, S(v) \, \langle v \rangle^q 
 + \int_{\ens{S}^{N-1}} 
     \left( \frac{\langle v \rangle^q - 
     \langle \bar v \rangle^q}{(\cos \theta/2)^{N+\gamma}} \right) 
     S(\bar{v}) \, b \, d\sigma ,
 \end{multline*}
which in turn yields a corresponding splitting of $I_{2}$ into three parts 
$I_{2,1} + I_{2,2} + I_{2,3}$. 

For the $I_{2,3}$ term, we use again the change of variable of cancellation lemmas backward and 
we use~(\ref{ll}) (for $p> N/(N+\gamma)$): 
 \begin{multline*}
 I_{2,3} \le C \, \int_{v,v_*} |D_*| S(v) \, |v-v_*|^\gamma \,  
               \left( \int_{\ens{S}^{N-1}} 
          \left| \langle v' \rangle^q - \langle v \rangle^q \right| 
              \, b \, d\sigma \right) \, dv \, dv_* \\
          \le C \, \left( \|S\|_{L^1_q} + \|S\|_{L^p} \right) 
                 \, \|D\|_{L^1_{q}}. 
 \end{multline*}

For the $I_{2,2}$ term,  using again 
 \[ \left| \frac{1}{(\cos \theta/2)^{N+\gamma}} -1 \right| \le C \, (1-\cos \theta), \]
 we get  
 \[  I_{2,2} \le C \,\left( \|S\|_{L^1_q} + \|S\|_{L^p} \right)  \, \|D\|_{L^1}. \] 

Finally, for the $I_{2,1}$ term, we assume first for simplicity $\nu < 1$ and we 
denote $\bar{v}_u = (1-u)\,v + u\,\bar{v}$. We Taylor-expand $S \langle \cdot \rangle^q$: 
 \[ I_{2,1} \le C \, \int_{\R^{2N} \times \ens{S}^{N-1} \times [0,1]}
       |\nabla(S \langle \cdot \rangle^q)|(\bar{v}_u) \, 
       |v-v_*|^{\gamma+1} \, (\tan \theta/2) \, b 
       \, |D_*| \, dv \, dv_* \, d\sigma \, du. \]
Hence we deduce 
 \begin{equation}
 I_{2,1} \le 
      C \, \left(\int_{\ens{S}^{N-1}}(\tan \theta/2) \, b \, d\sigma \right)
          \,  \|S\|_{W^{1,1}_{q+1 + \gamma}} 
           \, \|D\|_{L^1 _{q}}
 \end{equation}
when $\gamma +1 \ge 0$, and 
 \begin{equation}
 I_{2,1} \le 
      C \, \left(\int_{\ens{S}^{N-1}}(\tan \theta/2) \, b \, d\sigma \right)
          \, \left( \|S\|_{W^{1,1}_{q}} + \| \nabla S \|_{L^p} \right)   
           \, \|D\|_{L^1 _{q}}
 \end{equation}
with $p > N/(N+\gamma+1)$ else. The case $1 \le \nu <2$ can be treated in a similar 
way by Taylor-expanding at second order as in the previous proof. 



Combining all the previous estimates, we obtain 
 \[ \frac{d}{dt} \|D\|_{L^1 _q}\le C_s\, \|D\|_{L^1 _q}, \]
which concludes the proof of Theorem~\ref{theo:stab} and Proposition~\ref{propo:sta}. 

\section{Proof of the estimates on the propagation of smoothness} \label{eps}
\setcounter{equation}{0}

We now turn to the
\bigskip

{\bf{Proof of Theorem~\ref{theo:reg}}:}
We begin the proof (Subsection \ref{hpmsp}) under the assumptions {\bf H1}, {\bf H2} with $\nu <1$
and ({\bf H3-1}, {\bf H3-2} or {\bf H3-3}). We detail only the cases of assumptions 
{\bf H3-3} or {\bf H3-4}, since the two other cases are similar (and somewhat simpler).

\subsection{Hard potentials and mollified soft potentials} \label{hpmsp}

We assume that (for some $q\ge 2$), we have 
an initial datum $f(0,\cdot) \in W^{1,1}_q (\R^N)$, and we consider a solution 
$f$ to the spatially homogeneous Boltzmann equation~\eqref{el}.

We split $Q$ (and correspondingly $B$, $b$) like in the previous section into two operators
$Q_c$ and $Q_r$. 
Then, we compute (for $q\ge 2$), denoting by $\partial_v f$ any (first order) 
partial derivative of $f$ with respect to 
one of the components of $v$, the following quantity:
\begin{multline*}
 \frac{d}{dt} \int_{\R^N} |\partial_v f| \, \langle v \rangle^q \, dv 
 = \int_{\R^N}  Q_r(f,\partial_v f) \, \mbox{sgn}(\partial_v f) \, \langle v \rangle^q \, dv \\
 + \int_{\R^N}  Q_c(f,\partial_v f) \, \mbox{sgn}(\partial_v f) \, \langle v \rangle^q \, 
 dv  := I_1 + I_2 .
  \end{multline*}
We first consider the term $I_2$ corresponding to the cutoff part:
$$
 I_2 = \int_{v,v_*,\sigma} \Big[ f' _* \partial_v f' + f' \partial_v f' _* - f_* \partial_v f - f \partial_v f_* \Big] 
 \, \mbox{sgn}(\partial_v f) \, \langle v \rangle^q \, B_c  $$
 $$ \le \int_{v,v_*,\sigma} \Big[ f' _* |\partial_v f'| + f' |\partial_v f'| _* 
                  - f_* |\partial_v f| - f |\partial_v f_*| \Big]  \, \langle v \rangle^q \, B_c
 + 2\,\int_{v,v_*,\sigma} f |\partial_v f _*| \, \langle v \rangle^q \, B_c $$
 $$ = \int_{\R^N} Q_c(f,|\partial_v f|) \, \langle v \rangle^q \, dv 
 + 2\,\int_{v,v_*,\sigma} f |\partial_v f_*| \, \langle v \rangle^q \, B_c .    $$
Arguing as in the proof of the stability estimates, we get
  \begin{equation}\label{yyy}
  \frac{d}{dt} \int_{\R^N} |\partial_v f| \, \langle v \rangle^q \, dv \le C_\var \, \| \partial_v f\|_{L^1 _q} 
           - K\, \| \partial_v f\|_{L^1 _{q+\gamma_+}} 
	   \end{equation}
where $C_\var$ depends on  $\|f\|_{L^1_{q+\gamma}}$ and $\var$ (indeed, as 
explained in the proof of the stability estimates, since $\nu <1$ the constant $C_\var$ can be 
taken independent of $\var$), and $K$ only depends on a constant $C>0$  such that
  \[ \int_{\R^N} f \, |v-v_*|^\gamma \, dv_* \ge C \, \langle v \rangle^{\gamma} . \]

We now turn to the non cutoff part. We write
$$ I_1 =  \int_{v,v_*,\sigma} \Big[ f' _* \partial_v f' + f' \partial_v f' _* - f_* \partial_v f - f \partial_v f_* \Big] 
 \, \mbox{sgn}(\partial_v f) \, \langle v \rangle^q \, B_r  $$
 $$ = \int_{v,v_*,\sigma} \Big[ \mbox{sgn}(\partial_v f)' \, \langle v' \rangle^q 
     + \mbox{sgn}(\partial_v f)' _* \, \langle v' _* \rangle^q 
 - \mbox{sgn}(\partial_v f) \, \langle v \rangle^q - \mbox{sgn}(\partial_v f)_* \, \langle v_* \rangle^q \Big] 
  \, f_* \partial_v f  \,B_r $$
  $$ \le \int_{v,v_*,\sigma} \Big[ \mbox{sgn}(\partial_v f)(v'_*)  -  \mbox{sgn}(\partial_v f)(v_*) \Big] 
       \, f_* \partial_v f(v) \, \langle v_* \rangle^q \, B_r  $$
  $$ + \int_{v,v_*} \left| \int_\sigma \left[ \langle v' \rangle^q - \langle v \rangle^q \right] \, b_r \, d\sigma \right| \, 
        \left( f_* \, |\pa_v f(v)| + f \, |\pa_v f(v_*)| \right) \, |v-v_*|^\gamma =: I_{1,1} + I_{1,2}. 
    $$
The term $I_{1,2}$ is easily controlled thanks to Lemma~\ref{diffpoids}: 
$$
I_{1,2} \le \delta(\var) \, \| \pa_v f \|_{L^1 _{q+\gamma}} 
$$
with $\delta(\var) \to 0$ as $\var \to 0$. 
 
We now focus on the most difficult term $I_{1,1}$. 
Since the proof makes use, in the sequel, of an intricate kind of integration by parts, we write down
first the simple case when the dimension is $N=2$, for the sake of clarity. 
\par
In this case ($N=2$), 
we define the change of variables $v_* \mapsto w$ (for given $v$, $\theta$), where 
 $$ w = v'_* = \frac{v+v_*}2 - R_{\theta}\bigg( \frac{v-v_*}2 \bigg) $$
($R_\theta$ denotes the rotation of angle $\theta$), 
whose Jacobian determinant is clearly $(\cos \theta/2)^2$ and which can be inverted in
 $$ v_* = v_*(v,w,\theta) = \frac{R_{- \frac{\theta}2} w + \sin \frac{\theta}{2} \, R_{\frac{\pi}2} v }{\cos \frac{\theta}{2} } . $$
Using this change of variables in the first part of the integral (and the fact that 
$|v-w| = \cos \theta /2 \, |v-v_*|$), we obtain
 $$ 
 I_{1,1} \le \int_{v\in \R^2} \int_{w\in\R^2} \int_{\theta=0}^{\var} \mbox{sgn}(\partial_v f)(w) \,  
  \partial_v f (v)  \, b_r(\cos \theta) 
 $$
$$ 
\times \bigg\{   f\big( v_*(v,w,\theta) \big) \, \langle v_*(v,w,\theta) \rangle^q 
\, \big| v - v_*(v,w,\theta) \big|^{\gamma}  
\,\left(\cos \frac{\theta}{2}\right)^{-2}  $$
$$- f(w)\, \langle w \rangle^q \, |v-w|^{\gamma} \bigg\} \, d\theta \, dw \, dv  
$$
$$ 
\le  \int_{v\in \R^2} \int_{w\in\R^2} \int_{\theta=0}^{\var}  
\mbox{sgn}(\partial_v f)(w) \,  \partial_v f (v) \, b_r(\cos \theta) $$
$$\times \, \bigg[ \left(\cos \frac{\theta}{2}\right)^{-2-\gamma} - 1 \bigg] \,  f(w) 
     \langle w \rangle^q \, |v-w|^{\gamma} \, d\theta \, dw \, dv  
$$
$$ 
+  \int_{v\in \R^2} \int_{w\in\R^2} \int_{\theta=0}^{\var}
\mbox{sgn}(\partial_v f)(w) \,  \partial_v f (v) \, b_r(\cos \theta) \, 
\left(\cos \frac{\theta}{2}\right)^{-2-\gamma}\,
$$
$$  
\times \big\{  f\big( v_*(v,w,\theta) \big) -f(w) \big\} \, \langle w \rangle^q \, 
|v-w|^{\gamma} \, d\theta \, dw \, dv  
$$
$$
+ \int_{v\in \R^2} \int_{w\in\R^2} \int_{\theta=0}^{\var}
\mbox{sgn}(\partial_v f)(w) \,  \partial_v f (v) \, b_r(\cos \theta) \, 
\left(\cos \frac{\theta}{2}\right)^{-2-\gamma}\,
$$
$$  
\times \big\{  \langle v_*(v,w,\theta) \rangle^q -\langle w \rangle^q \big\} \, f(v_*(v,w,\theta) ) \, 
|v-w|^{\gamma} \, d\theta \, dw \, dv  =: I_{1,1,1} + I_{1,1,2} + I_{1,1,3}. 
$$
The first term $I_{1,1,1}$ is controlled thanks to 
$$
\bigg| \left(\cos \frac{\theta}{2}\right)^{-2-\gamma} - 1 \bigg| \le C \, (1 - \cos \theta)
$$
which yields 
$$
I_{1,1,1} \le \mbox{{\rm cst}} \, \| \pa_v f \|_{L^1 _q}. 
$$
The third term $I_{1,1,3}$ is controlled thanks to the argument of Lemma~\ref{diffpoids}: 
$$
I_{1,1,3} \le \delta(\var) \, \| \pa_v f \|_{L^1 _{q+ \gamma}}. 
$$
Finally we use integration by parts (according to $v$) on the second term $I_{1,1,2}$ (and the fact 
that the differential in $v$ of $v_*(v,w,\theta)$ has bound $\tan \theta /2$):
$$ 
I_{1,1,2} \le  \int_{v\in \R^2} \int_{w\in\R^2} \int_{\theta=0}^{\var}
 b_r(\cos \theta) \,  \left(\cos \frac{\theta}{2} \right)^{-3-\gamma} \,  \left(\sin \frac{\theta}{2}\right) \,
f(v) 
$$
$$ 
\times \, \left| \partial_v f  \right| \big(  v_*(v,w,\theta) \big) \, |v-w|^{\gamma} 
\, \langle w \rangle^q \, d\theta \, dw \, dv  
$$ 
$$ 
+ \gamma \, \int_{v\in \R^2} \int_{w\in\R^2} \int_{\theta=0}^{\var}
b_r(\cos \theta) \,  \left(\cos \frac{\theta}{2}\right)^{-2-\gamma} \, f(v)
$$
$$ 
\times \,\left|  f\big( v_*(v,w,\theta) \big) -f(w) \right| \, |v-w|^{\gamma-1} \,  \langle w \rangle^q 
  \, d\theta \, dw \, dv  
$$
$$ 
\le  \int_{v\in \R^2} \int_{v_*\in\R^2} \int_{\theta=0}^{\var} b_r(\cos \theta) \, \left( \tan \frac{\theta}{2} \right) \,
f(v)\, \langle v' _* \rangle^q  \, |v-v_*|^{\gamma} \, |\partial_v f(v_*)|  \, d\theta \, dv_* \, dv  
$$
$$ 
+ \gamma \, \int_{v\in \R^2} \int_{v_*\in\R^2} \int_{\theta=0}^{\var}  \int_{u=0}^1
b_r(\cos \theta) \,  \left(\cos \frac{\theta}{2}\right)^{-3-\gamma} \,  \left( \sin \frac{\theta}{2} \right) 
\, f(v)\,  |v-w|^\gamma  \, \langle w \rangle^q 
$$
$$ 
\times \, \left| \partial_v f \right| \big( (1-u)\,w + u\, v_*(v,w,\theta) \big) \, du \, d\theta \, dw \, dv
$$
$$ 
\le \delta(\var)\,  \| \partial_v f\|_{L^1_{q+\gamma}}, 
$$
where $\lim_{\var \to 0} \delta(\var) =0$. Here, we have assumed 
 that $ \| f \|_{L^1_{q+\gamma}}$ is bounded (this norm is known to be 
at least propagated, Cf. \cite{Desv:93} for example).  
\par
Then, using estimate (\ref{yyy}), we see  that (optimizing $\var$), the proof of the theorem is complete. 

Let us now explain how to deal with the general case of dimension $N \ge 2$.

In the formula for $I_{1,1}$ 
\[ I_{1,1} = \int_{v,v_*,\sigma} \Big[ \mbox{sgn}(\partial_v f)(v'_*)  -  \mbox{sgn}(\partial_v f)(v_*) \Big] 
       \, f_* \partial_v f(v) \, \langle v_* \rangle^q \, B_r,  \]
we use the change of variable $v_* \mapsto w$ (for given $v,\sigma$) where
\[ w = v'_* = \frac{v+v_*}2 - \frac{|v-v_*|}2\,\sigma, \]
which can be inverted in
\[ v_* = v_*(v,w,\sigma) = 2\,w - v + \frac{|v-w|^2}{(v-w)\cdot\sigma}\,\sigma, \]
and whose Jacobian determinant is
\[ J = 2^{N-1}\, \frac{|v-w|^2}{((v-w)\cdot\sigma)^2}. \]
We get 
\[ I_{1,1} \le \bigg| \int_{v,w,\sigma}{\mbox{sgn}}(\pa_v f)(w)\, 
     f\left(v_*(v,w,\sigma)\right) \, 
      \pa_v f(v) \, \bigg( \frac{|v-w|^2}{(v-w)\cdot\sigma} 
       \bigg)^{\gamma} \, |v-w|^\gamma \]
\[ \times\, \langle v_*(v,w,\sigma) \rangle^q \, b_r \left(-1 + 2\,\left(\frac{v-w}{|v-w|}\cdot\sigma \right)^2 \right) \, 
   \frac{2^{N-1}}{\left(\frac{v-w}{|v-w|}\cdot\sigma\right)^2} \, d\sigma \, dw \, dv \]
\[ \qquad \qquad - \int_{v,w,\sigma} {\mbox{sgn}}(\pa_v f)(w)\, f(w)\, \pa_v f(v) \, |v-w|^{\gamma}
\, \langle w \rangle^q \, b_r \left(\frac{v-w}{|v-w|}\cdot\sigma \right)\, d\sigma \, dw \, dv \bigg|, \]
and using the decomposition of the unit vector $\sigma$:
\[ \sigma = \cos\theta \, \frac{v-w}{|v-w|} + \sin\theta \, n, \]
where $ n \in S^{N-1} \cap (v-w)^{\bot}$ and $0\le \theta \le \var$ (we look here to the non cutoff part of the cross-section), we end up with
\[ I_{1,1} \le \bigg| \int_{v,w} \int_{\theta = 0}^{\var} 
 \int_{n\in S^{N-1} \cap (v-w)^{\bot}} {\mbox{sgn}}(\pa_v f)(w)\, f\big(w+\tan\theta\,|v-w|\,n\big) \, \pa_v f(v) \]
\[\times \, |v-w|^{\gamma}\,
 \big\langle w+\tan\theta\,|v-w|\,n \big\rangle^q \, \frac{b_r\big(\cos (2\theta)\big)}{\cos ^{2+\gamma} \theta} \, 
\, 2^{N-1}\,\sin^{N-2} \theta \, dn \, d\theta \, dw \, dv \]
\[ - \int_{v,w} \int_{\theta = 0}^{\var}
 \int_{n\in S^{N-1} \cap (v-w)^{\bot}} {\mbox{sgn}}(\pa_v f)(w)\, f(w)\, 
 \pa_v f(v)\, |v-w|^{\gamma} \]
 \[ \hspace{6cm}  
  \times\,  \langle w \rangle^q \, b_r(\cos\theta) \, \sin^{N-2} \theta \, dn \, d\theta \, dw \, dv \bigg|. \]
After integration by part (in $v$), and the use of the change of variables $\delta = 2\theta$ in 
the first term, we obtain
\[ I_{1,1} \le \bigg| \int_{v,w} \int_{\delta=0}^{\var} 
\int_{n\in S^{N-1} \cap (v-w)^{\bot}}
{\mbox{sgn}}(\pa_v f)(w)\, |v-w|^{\gamma} \] 
 \[ \bigg[ f\big(w+\tan \delta/2 |v-w|\,n\big)\, \left(\cos \frac{\delta}2\right)^{-\gamma - N} \, 
    \big\langle w+\tan\delta/2 \,|v-w|\,n \big\rangle^q -f(w) \, \langle w \rangle^q \bigg] \] 
 \[  \pa_v f (v)\, b_r(\cos \delta) \, \sin^{N-2} \delta \, dn \, d\delta \, dw \, dv \bigg| \]
\[ +  \bigg| \int_{v,w} \int_{\delta=\var}^{2\var} \int_{n\in S^{N-1} \cap (v-w)^{\bot}} 
{\mbox{sgn}}(\pa_v f)(w)\, |v-w|^{\gamma}\,
 \big\langle w+\tan\delta/2 \,|v-w|\,n \big\rangle^q \]
\[ \times \, f\big(w+\tan \delta/2 |v-w|\,n\big) \, (\pa_v f)(v)\, \left(\cos \frac{\delta}2\right)^{-\gamma - N} \, 
   b_r(\cos \delta) \, \sin^{N-2} \delta \, dn \, d\delta \, dw \, dv \bigg| \] 
\[  =: I_{1,1,1} + I_{1,1,2} + I_{1,1,3}. \] 
Then for the first term $I_{1,1,1}$ we have 
\[ I_{1,1,1} \le \gamma \, \bigg| \int_{v,w} \int_{\delta=0}^{\var}
\int_{n\in S^{N-1} \cap (v-w)^{\bot}}
f(v)\, {\mbox{sgn}}(\pa_v f)(w)\, |v-w|^{\gamma-1} \] 
 \[ \bigg[ f\big(w+\tan \delta/2 |v-w|\,n\big)\, \left(\cos \frac{\delta}2\right)^{-\gamma - N} \, 
    \big\langle w+\tan\delta/2 \,|v-w|\,n \big\rangle^q -f(w) \, \langle w \rangle^q \bigg] \] 
 \[  b_r(\cos \delta) \, \sin^{N-2} \delta \, dn \, d\delta \, dw \, dv \bigg| \]
from which we deduce straightforwardly as before 
 \[ I_{1,1,1} \le \delta(\var) \, \| \pa_v f \|_{L^1 _{q+ \gamma}} \]
by using the control on $(\cos^{-\gamma - N} \delta/2 -1)$ and 
Taylor-expanding $f$ and $\langle \cdot \rangle^q$. 
 
For the second term $I_{1,1,2}$ we have 
 \[ I_{1,1,2} \le \bigg| \int_{v,w} \int_{\delta=0}^{\var} 
 f(v)\, {\mbox{sgn}}(\pa_v f)(w)\, |v-w|^{\gamma} \]
 \[ \hspace{3cm} \times \pa_v \bigg( \int_{n\in S^{N-1} \cap (v-w)^{\bot}} 
\big(f \langle \cdot \rangle^q \big)\big(w+\tan \delta/2 |v-w|\,n\big)\, dn \bigg) \]
\[ \hspace{6cm} \times\, \left(\cos \frac{\delta}2\right)^{-\gamma - N} \, b_r(\cos \delta)
 \, \sin^{N-2} \delta \, d\delta \, dw \, dv \bigg| \]
 \[ \le \mbox{{\rm cst}} \, \bigg| \int_{v,w} \int_{\delta=0}^{\var} 
 f(v)\, {\mbox{sgn}}(\pa_v f)(w)\, |v-w|^{\gamma} \]
 \[ \hspace{3cm} \times \bigg( \int_{n\in S^{N-1} \cap (v-w)^{\bot}} 
\big| \pa_v \big(f \langle \cdot \rangle^q \big) \big| \big(w+\tan \delta/2 |v-w|\,n\big)\, dn \bigg) \]
\[ \hspace{6cm} \times\, \left(\cos \frac{\delta}2\right)^{-\gamma - N} \, \tan \delta /2 \, b_r(\cos \delta)
 \, \sin^{N-2} \delta \, d\delta \, dw \, dv \bigg| \]
(the derivative taken on a sphere depending on $v$ has been treated by taking
local coordinates) from which we deduce 
$$
I_{1,1,2} \le \delta(\var) \, \| \pa_v f \|_{L^1 _{q+ \gamma}}. 
$$

For the third term we have 
\[ I_{1,1,3} \le  \bigg| \int_{v,w} \int_{\delta=\var}^{2\var} \int_{n\in S^{N-1} \cap (v-w)^{\bot}} \]
\[  \pa_v f(v)\, {\mbox{sgn}}(\pa_v f)(w)\, |v-w|^{\gamma}\,
 \big\langle w+\tan\delta/2 |v-w|\,n \big\rangle^q \,  f\big(w+\tan\delta/2 |v-w|\,n \big) \] 
 \[ \hspace{3cm} \times \, \left(\cos \frac{\delta}2\right)^{-\gamma - 2} \, b_r(\cos \delta)
 \, \sin^{N-2} \delta \, dn \, d\delta \, dw \, dv \bigg| \]
from which we deduce thanks to the cutoff (by coming back to the classical variables)
$$
I_{1,1,3} \le C_\var \, \| \pa_v f \|_{L^1 _\gamma} 
$$
(for some constant $C_\var$ blowing-up as $\var \to 0$).

These estimates (where we have used the boundedness of weighted $L^1$ norms 
like in dimension $2$)  enables to complete the proof of Theorem~\ref{theo:reg} 
(for hard and smoothed soft potentials) also in dimension bigger than 2.

\subsection{Non mollified soft potentials}
\par
We now prove Theorem~\ref{theo:reg} under assumptions {\bf H1}, {\bf H2} with $\nu <1$
and {\bf H3-4}.  
The additional difficulty here is the singularity of the collision kernel for 
small relative velocity. As pointed out in the proofs of stability estimates, 
this suggests to use some $L^p$ norm with $p> N /(N+\gamma)$ of $f$,
in order to control some convolution terms of the form 
$$
v \mapsto \int_{|v-v_*| \le 1} |v-v_*|^\gamma \, f(v_*) \, dv_*. 
$$

We first state a result showing that we are able to propagate $L^1$ norms of the gradient 
as soon as we have some time integrability of some $L^p$ and $L^1$ moments. 
\begin{proposition} \label{proppotmous3} 
Let $B$ be a collision kernel satisfying assumptions {\bf H1}, {\bf H2} with $\nu <1$ 
and {\bf H3-4}. Let $f$ be a solution on $[0,T]$
to the corresponding spatially homogeneous Boltzmann 
equation (\ref{el}). Suppose that for $p>N/(N+\gamma)$ and $q\ge 2$, one has
 $f\in L^1([0,T]; L^p \cap L^1_q (\R^N))$.
Assume also that  $f(0,\cdot) \in W^{1,1}_q(\R^N)$. Then
$f \in L^{\infty}([0,T]; W^{1,1}_q(\R^N))$.
\end{proposition}

\smallskip\noindent{\sl Proof of proposition \ref{proppotmous3}.} 
We compute (for $q \ge 2$)
\begin{multline*} 
{\frac{d}{dt}} \int |\partial_v f|\,  \langle v \rangle^q  \,dv  
=  \int_{v,v_*,\sigma} \big[ f' _* \partial_v f' + f' \partial_v f' _* - f_* \partial_v f - f \partial_v f_* \big] 
 \, \mbox{sgn}(\partial_v f) \, \langle v \rangle^q \, B  \\
 \le \int_{v,v_*,\sigma} \big[  \mbox{sgn}(\partial_v f)(v'_*)  -  \mbox{sgn}(\partial_v f)(v_*) \big]\, f_*
     \partial_v f \, \langle v_* \rangle^q \, B \\
     + \int_{v,v_*} \big( 
      f_* \, |\partial_v f (v)| + f \, | \pa_v f (v_*)| \big) 
      \, \left| \int_\sigma \left[ \langle v' \rangle^q - \langle v \rangle^q \right] \right| \, B 
   \qquad   =: I_1 + I_2.
\end{multline*}
The term $I_2$ is immediately controlled by Lemma~\ref{diffpoids}: 
$$ 
I_2 \le C \, \| \partial_v f \|_{L^1_q}\, (\| f \|_{L^1} + \| f \|_{L^p})
$$
for some constant $C >0$. Then we focus on the term $I_1$ and we split it as $I_1 ^c + I_1 ^r$ 
according to the decomposition of the collision kernel $B= B_c + B_r$ as before. 

We first note that (for some constant $C_\var$ possibly blowing-up when $\var \to 0$)
$$ I_1 ^c \le \bigg| \int_{v,v_*,\sigma} \big[  \mbox{sgn}(\partial_v f)(v'_*)  -  \mbox{sgn}(\partial_v f)(v_*) \big]\, f_*
     \partial_v f \, \langle v_* \rangle^q \, B_c \bigg| $$
$$   \le C_{\var}\, \| \partial_v f \|_{L^1_q}\, (\| f \|_{L^1} + \| f \|_{L^p}) . $$

Then, we only detail the treatment of the non cutoff case in the simple case when the dimension is $N=2$
(the case of dimension $N>2$ can be treated like in Subsection~\ref{hpmsp}, 
where hard potentials are considered).
Using the change of variables $v_* \mapsto w$ (for a given $v$, $\theta$), we end up with 
\begin{multline*} 
I_1 ^r   \le  \bigg| \int_{v\in \R^2} \int_{w\in\R^2} \int_{\theta=0}^{\var}  
	       \mbox{sgn}(\partial_v f)(w) \,  \partial_v f (v)  \, \langle w \rangle^q \\
	       \times\,  b_r(\cos \theta) \, \bigg( \left(\cos \frac{\theta}{2}\right)^{-2-\gamma} - 1 \bigg) 
	       \,  f(w)\, |v-w|^{\gamma} \, d\theta \, dw \, dv  \bigg| \\
	       +  \bigg| \int_{v\in \R^2} \int_{w\in\R^2} \int_{\theta=0}^{\var}
	       \mbox{sgn}(\partial_v f)(w) \,  \partial_v f (v)   \, b_r(\cos \theta) \, 
	       \left(\cos \frac{\theta}{2}\right)^{-2-\gamma} \\ 
	       \times \big\{  f\big(v_*(v,w,\theta) \big)  \, \langle v_*(v,w,\theta) \rangle^q 
	       - f(w)  \, \langle w \rangle^q\big\} \, |v-w|^{\gamma}  \, d\theta \, dw \, dv \bigg| \\ 
	       \hspace{-3cm} \le  \, \mbox{{\rm cst}} \,   \int_{v\in \R^2} \int_{w\in\R^2} 
	        | \partial_v f (v)| \, \langle w \rangle^q 
	        \,  f(w)\, |v-w|^{\gamma} \, dw \, dv \\ 
	       \hspace{-3cm} + \int_{v\in \R^2} \int_{w\in\R^2} \int_{\theta=0}^{\var}
	        b_r(\cos \theta) \,  \left(\cos \frac{\theta}{2}\right)^{-2-\gamma} \,  \sin \frac{\theta}{2} \\ 
		\times \, f(v)\, | \partial_v (f \, \langle \cdot \rangle^q )|\big( v_*(v,w,\theta) \big) \,
		|v-w|^{\gamma}  \,d\theta \, dw \, dv  \\
                \hspace{-3cm} + |\gamma| \, \int_{v\in \R^2} \int_{w\in\R^2} \int_{\theta=0}^{\var}
           b_r(\cos \theta) \,  \left(\cos \frac{\theta}{2}\right)^{-2-\gamma} \, f(v) \\ 
     \times \,\big| f\big(v_*(v,w,\theta) \big)  \, \langle v_*(v,w,\theta) \rangle^q 
	       - f(w)  \, \langle w \rangle^q \big| \, 
	       |v-w|^{\gamma-1} \,d\theta \, dw \, dv,  
\end{multline*}
and we deduce (using Taylor expansion as before on the last term)  
$$
I_1 ^r     \le \, \mbox{{\rm cst}} \, \| \pa _v f \|_{L^1 _q} \, \Big(  \|f\|_{L^1 _{q}} +  \|f\|_{L^p} \Big) 
$$
which concludes the proof of the proposition by some Gronwall argument. 
\qed

\bigskip

Then, the $L^p$ norms can easily be shown to be locally propagated in time by direct energy methods:

\begin{proposition} \label{proppotmous} 
Let $B$ be a collision kernel satisfying assumptions {\bf H1}, {\bf H2} with $\nu <1$ 
and {\bf H3-4}. Let $f$ be a solution on $[0,T]$
to the corresponding spatially homogeneous Boltzmann 
equation (\ref{el}).  
Suppose that $f(0,\cdot) \in L^p(\R^N)$ for some $p\in ]N/(N+\gamma), +\infty]$.
Then, there exists a time $T^* \in ]0,T]$
explicitly computable such that
$f \in L^{\infty}([0,T^*]; L^{p}(\R^N))$.
\end{proposition}


\smallskip\noindent{\sl Proof of proposition \ref{proppotmous}.} 
We compute for some $N/(N+\gamma)< p < +\infty$ 
(and with $C$ denoting a constant which does not depend on $p$)
$$ {\frac{d}{dt}} \int_{v\in\R^N} |f(t,v)|^p\, dv = p\, \int_{v,v_*,\sigma} 
\bigg(f^{p-1}\,f'\,f^{'*} - f^p \,f^* \bigg) \, B $$
$$ = p\, \int_{v,v_*,\sigma} \bigg((f')^{p-1}\,f\,f^{*} - f^p \,f^* \bigg) \, B $$
$$ \le  p\, \int_{v,v_*,\sigma} \bigg( \frac1p\, f^p + (1 - \frac1p)(f')^{p} - f^p\bigg) \, f^* \, B $$
$$ \le (p-1)\, \int_{v,v_*,\sigma} \bigg(  (f')^{p} - f^p\bigg) \, f^* \, B $$
$$ \le C\,(p-1)\,  \int_{\theta=0}^{\pi/2} \left( \left( \cos \frac{\theta}2\right)^{-2 + \gamma} - 1\right) 
\, b(\cos \theta) \,d\theta \,\, 
\int_{v,v_*} f^p \, f^*  \, |v - v_*|^{\gamma}  $$
$$ \le  C\, (p-1)\, \int_{v}  f^p \, \bigg( \int_{|v_* - v| \le 1}   f^*  \, |v - v_*|^{\gamma}  \, dv_*
+  \int_{|v_* - v| \ge 1}   f^*  \, |v - v_*|^{\gamma} \,dv_*\bigg)  \, dv $$
$$ \le  C\, (p-1)\, \int_{v}  f^p \, \bigg[ C\, \bigg( \int (f^*)^p \bigg)^{1/p}  + \int f^* \bigg] $$
$$ \le  C\, (p-1)\, \bigg[ \bigg( \int f^p \bigg)^{ 1 + 1/p} +  \int f^p \bigg] $$
(remember in theses computations that the collision kernel $B$ is taken in symmetrized form 
$B_{sym}$ with support in $[0,\pi/2]$). As a consequence, 
$$ 
\frac{d}{dt} \| f \|_{L^p} \le C\, \left(1 - \frac1p\right)\, \Big( \| f \|_{L^p} + \| f \|_{L^p}^2 \Big), 
$$
and we get also (passing to the limit when $p \to + \infty$),
$$ 
\frac{d}{dt} \| f \|_{L^{\infty}}  \le C \,  \Big( \| f\|_{L^{\infty}} + \| f \|_{L^{\infty}}^2\Big) . 
$$
This ends the proof of Proposition~\ref{proppotmous}.
\qed
\bigskip

Finally, the $L^1$ moments can be shown to be locally propagated, using the previous 
$L^p$ estimates in case of strong singularity at small relative velocity: 

\begin{proposition} \label{proppotmous2}
Let $B$ be a collision kernel satisfying assumptions {\bf H1}, {\bf H2} with $\nu <1$ 
and {\bf H3-4}. Let $f$ be a solution on $[0,T]$
to the corresponding spatially homogeneous Boltzmann 
equation (\ref{el}). 

Suppose first that $\gamma + 1 \ge 0$. 
Then for any $q \ge 2$, 
if $f(0,\cdot) \in L^1 _q(\R^N)$, one has $f\in L^{\infty}([0,T];L^1 _q (\R^N))$.

Suppose then that $\gamma + 2 \ge 0$. 
Then for any $q \ge 4$, 
if $f(0,\cdot) \in L^1 _q(\R^N)$, one has $f\in L^{\infty}([0,T];L^1 _q (\R^N))$.

Suppose finally that $\gamma +2 < 0$ and  
$f(0,\cdot) \in L^p(\R^N) \cap L^1_q(\R^N)$ for some 
$N/(N+\gamma+2) < p \le +\infty$, $q \ge 4$. 
Then, there exists a time $T^* \in ]0,T]$ explicitly computable such that
$f\in  L^{\infty}([0,T^*]; L^{p}(\R^N)\cap L^1_q(\R^N))$.
\end{proposition}

\smallskip\noindent{\sl Proof of proposition \ref{proppotmous2}.} 
We compute (for $q\ge 2$)
$$ {\frac{d}{dt}} \int_{v\in\R^N} f(t,v) \, \langle v \rangle^q \, dv
= \int_{v,v_*,\sigma} \langle v \rangle^q ( f'\,f^{'*} - f \,f^*)\,  B $$
$$ \le   \int_{v,v_*}  \bigg[ \int_{\sigma} \Big[ 
\langle v' \rangle^q - \langle v \rangle^q \Big] \, b \, d\sigma \bigg]
\, f \, f^*  \, |v - v_*|^{\gamma} . $$
Using Lemma~\ref{diffpoids}, we get for $q \ge 2$:
 $$ {\frac{d}{dt}} \int_{v\in\R^N} f(t,v) \, \langle v \rangle^q \, dv 
\le \mbox{{\rm cst}} \, \int_{v,v_*} \Big[ \langle v_* \rangle^{q-1} + \langle v \rangle^{q-1} \Big] \, f \, f^* 
\, |v - v_*|^{\gamma+1} $$
and for $q \ge 4$: 
$$ {\frac{d}{dt}} \int_{v\in\R^N} f(t,v) \, \langle v \rangle^q \, dv 
\le \mbox{{\rm cst}} \, \int_{v,v_*} \Big[ \langle v_* \rangle^{q-2} + \langle v \rangle^{q-2} \Big] \, f \, f^* 
\, |v - v_*|^{\gamma+2}. $$
Hence this is controlled by 
$$
\mbox{{\rm cst}} \, \| f \|_{L^1_q} \, \| f \|_{L^1}
$$
if $q \ge 2$ and $\gamma + 1 \ge 0$, or $q \ge 4$ and $\gamma + 2 \ge 0$ 
(which concludes the proof of Proposition~\ref{proppotmous2} in those cases
 immediately by some Gronwall argument). It is also controlled by 
$$ \mbox{{\rm cst}} \, \| f \|_{L^1_q} \, \bigg( \| f \|_{L^1} + \| f \|_{L^{p}} \bigg) $$
with $p > N/(N+\gamma +2)$ if $q\ge 4$ and $\gamma +2 < 0$. 
This concludes also the proof of Proposition~\ref{proppotmous2} in this case:  
using Proposition \ref{proppotmous}, we see that the $L^{1}_q(\R^N)$ norm of $f(\cdot)$
is uniformly bounded on $[0,T^*]$, where $T^*$ is the same as in 
Proposition~\ref{proppotmous}. 
\qed 
\bigskip

In order to conclude the proof of Theorem~\ref{theo:reg}, it remains to prove that the $L^p$ estimate
is global in time in the case when $\gamma \in (-\nu, 0]$.
\par
In order to do so, we shall use the regularizing effect of non cutoff interactions. 
The entropy {\it a priori} estimate ensures that the entropy production associated 
to the solution $f$ of equation (\ref{el})  is bounded in $L^1([0,T])$ 
(uniformly with respect to $T \in [0, +\infty)$). This means that
\[ \int_0^T \int_{v,v_*,\sigma} f_*\, f\, \log(f/f') \, B <+\infty . \]
  Since the entropy production is monotonous (increasing) with respect
to the cross-section $B$, we  see that
\[ \int_0^T \int_{v,v_*,\sigma} f_*\, f\, \log(f/f') \,\max(B,1) <+\infty . \]
Using the results in~\cite{ADVW}, 
we get 
$$
\int_0 ^T \big\| \sqrt{f} \big\|^2_{H^{\nu/2}(|v| \le R)} \, dt \le C_R \, (1+ \| f \|_{L^1}^2 \, T).
$$
By tracking the constant in the proof of this estimate, one finds $C_R = \mbox{{\rm cst}} \, R^{|\gamma|}$. 
Using a Sobolev embedding (remember that the Sobolev constant does not depend on $R$),
 one deduces 
$$
\int_0 ^T \big\| \sqrt{f} \big\|^2_{L^{2N/(N-\nu)}(|v| \le R)} \, dt \le \mbox{{\rm cst}} \, R^{|\gamma|} \, (1+T) ,
$$
and finally
$$
\int_0 ^T \| f \|_{L^{N /(N-\nu)}(|v| \le R)} \, dt \le \mbox{{\rm cst}} \,R^{|\gamma|} \, (1+T).  
$$
We now use a dyadic decomposition of the velocity space. For any $k \ge 0$, we have 
$$
\int_0 ^T \| f \|_{L^{N /(N-\nu)}(2^k \le |v| \le 2^{k+1})} \, dt \le C\, 2^{k |\gamma|} \, (1+T).  
$$
Therefore, for any $\alpha \in \R$, 
$$
\int_0 ^T \| f \|_{L^{N /(N-\nu)} _{\gamma-\alpha} (2^k \le |v| \le 2^{k+1})} \, dt \le 
C \, 2^{-\alpha k} \, (1+T).  
$$
By summing this estimate for $k=0, \dots, +\infty$, we get when $\alpha>0$:
\begin{equation}\label{lpe}
\int_0 ^T \| f \|_{L^{N/(N-\nu)} _{\gamma-\alpha} (\R^N)} \, dt \le C \, (1+T).
\end{equation}

We now use the fact that $\gamma > -\nu$. Since in particular $\gamma > -2$, we can use 
Proposition~\ref{proppotmous2} in order to propagate the $L^1 _q$ moments for any $q \ge 4$, and get
\begin{equation}\label{lunq}
\int_0 ^T \| f \|_{L^1_q (\R^N)} \, dt \le C \, (1+T).
\end{equation}

By interpolation between estimates (\ref{lpe}) and (\ref{lunq}), we see that
\begin{equation}\label{lunqq}
\int_0 ^T \| f \|_{L^p (\R^N)} \, dt \le C \, (1+T)
\end{equation}  for any $1 < p < N/(N-\nu)$ (if 
$f(0,\cdot) \in L^1 _q$ with $q$ big enough).
But $N/(N+\gamma) < N/(N-\nu)$ (since $\gamma > -\nu$), so that we can take
$p > N/(N+\gamma)$, and the assumptions of Proposition~\ref{proppotmous3} hold for all times. Finally, the required smoothness
is global in time.

\Remark Note that the threshold $\gamma = -\nu$ in this proof does not seem to be only a 
coincidence. Indeed, as explained in the introduction, for inverse power-laws interaction potentials in 
dimension $3$, it corresponds to moderately soft potentials, that is potentials $V(r)=r^{-s}$ with 
$s > 2$. This critical value also corresponds to the threshold below which there is no spectral gap for 
the linearized operator (and therefore no expected exponential relaxation rate towards equilibrium), 
below which it is not known how to show slowly growing bounds on the $L^1$ moments (that is growing 
more slowly than exponentially), below which it is not known how to build global smooth solutions (as pointed 
out in this paper). Note that Fournier in~\cite{Fournier} also has the same threshold. 

An interesting issue to be explored is to relate this focus point of so many mathematical difficulties to the 
physical considerations: the threshold $s=2$ is called ``Manev interaction" in the physical literature
(see for instance~\cite{BDIVI,BDIVII}), and a dimensional analysis in~\cite{BDIVI} shows that for $s>2$,
the Boltzmann collision term prevail on the mean-field term, whereas for $1< s < 2$, the Boltzmann 
collision term should be negligible in front of the mean-field term.

\end{document}